\documentstyle[12pt]{article}
\newtheorem{lemma}{Lemma}
\newtheorem{theorem}{Theorem}
\newtheorem{proposition}{Proposition}
\newcommand{\QED}{{\hfill$\Box$\medskip}}
\newcommand{\k}{K\"{a}hler }
\newcommand{\ke}{K\"{a}hler-Einstein }
\newtheorem{definition}{Definition}
\newtheorem{corollary}{Corollary}

\begin{document}

\begin{center}
{\Large \bf Indefinite K\"{a}hler-Einstein Metrics on Compact 
Complex Surfaces}\\
\vspace{1 cm}
Jimmy Petean\\
Department of Mathematics, SUNY\\
Stony Brook, NY 11794-3651\\
e-mail: jimmy@math.sunysb.edu
\end{center}
\vspace{1 cm}
Mathematics Subject Classification (1991) 53C50, 32J15, 53C55
\vspace{1 cm}
\begin{abstract}Indefinite \k  solutions
of the Einstein equations are studied, and it is almost completely
determined which compact complex surfaces admit such metrics.
\end{abstract}

\section{Introduction}
 
A pseudo-Riemannian metric on a smooth manifold is called Einstein
if the Ricci tensor of the Levi-Civita connection equals a scalar
multiple of the metric.
This equation first appeared as a particular case (vacuum) of the
Einstein field equation `with cosmological constant'. 
This field equation was introduced by Einstein in order to describe the 
influence of matter and electromagnetism on spacetime; 
solutions of the Einstein equations would then describe
`vacuum spacetimes'.
In this physical context, one naturally considers Lorentzian
type metrics. But mathematicians became mainly interested in the
Riemannian case, and positive-definite solutions of the Einstein
equations have been intensively studied in the last decades.
While we still do not know much about general Riemannian solutions,
strong results have been obtained about the existence of
(positive definite) K\"{a}hler-Einstein metrics (see \cite{aubin},
\cite{yau}). In \cite{Besse} the reader
can find a detailed  discussion of these topics and an extensive
list of references.

In this paper we will consider indefinite \ke metrics. Our main result 
will be a
classification of the compact complex surfaces which admit such
metrics.

Let us begin by considering a compact complex manifold $(M^{2n} ,J)$. 
Here $M$ is a $2n$-dimensional
smooth compact manifold and $J$ is an integrable almost complex structure
on $M$. If $n=2$, $M$ is called a (compact) complex surface.

A pseudo-Riemannian metric $g$ on $M^{2n}$ is said to be Hermitian 
(or $J$-compatible) if $g(x,y)=g(Jx,Jy)$ for all $x,y$. At any
point of the manifold one can choose an orthogonal basis of the tangent space
of the form $\{ x_1 ,Jx_1 ,...,x_n ,Jx_n \} $; so, if $g$ is Hermitian,
its signature is of the form $(2k,2l)$. 
In particular if $M$ is a complex surface,
any indefinite Hermitian metric on $M$ has signature $(2,2)$.

If $g$ is a Hermitian pseudo-Riemannian metric then $\omega (x,y)=g(Jx,y)$ is
a 2-form, called the \k form of $g$.

\begin{definition} A Hermitian pseudo-Riemannian metric $g$ 
is called \k if its \k form
is closed. In particular, if $g$ is not positive or negative definite,
it is called an indefinite \k  metric.
\end{definition}

Consider now the Levi-Civita connection $\nabla $
of $g$ on $M$. 
Assume that $g$ is K\"{a}hler; then $J$ is parallel with respect to 
$\nabla $. This is usually
stated only in the Riemannian case, but it is not difficult to check that
it is also valid in the indefinite case (by exactly  the same proof).
Let $Ric$ be the Ricci tensor of $\nabla $. 
Then $Ric$ is $J$-invariant and hence 
$\rho (x,y)=Ric(Jx,y)$ is a 2-form. It is called the Ricci form of $g$.
It is also true in the indefinite case that $-i\rho $ is the curvature
of the canonical line bundle of $M$ (the bundle of holomorphic 2-forms);
the proof is the same as in the Riemannian case. In particular $\rho $
is closed and the de Rham class $[\rho /2\pi  ]$ is equal to the first 
Chern class of $M$ in cohomology with real coefficients.

\begin{definition} An indefinite \k  metric $g$ on $M$ is called 
indefinite \ke if there exists $\lambda \in {\bf R}$ such that
$Ric=\lambda g$ (or $\rho =\lambda \omega $). In this case $\lambda $
is called the Einstein constant.
\end{definition}

If $g$ is an indefinite \ke metric on $M$ and $k\in {\bf R}$,
 then $\hat{g} =kg$ is also an indefinite \ke metric (even if $k<0$).
The \k form of $\hat{g}$ is $\hat{\omega }=k\omega $ while the
Ricci form is $\hat{\rho }=\rho $. If $\omega = \lambda \rho$,
then $\hat{\omega} = \lambda k\hat{\rho}$. Without loss of
generality, we may therefore assume that $\lambda $ is either 0 or 1.

\ \ 

Indefinite \ke metrics on compact complex surfaces is the object
of study of this paper. Let us begin by constructing the simplest
examples.  

\ \

{\em Complex Tori}: Let $M= {\bf C}^2 / \Lambda $ be a complex 2-dimensional
torus. Let $z_1 ,z_2 $ be the standard coordinates on ${\bf C}^2 $. The
1-forms
$dz_1 ,dz_2 ,d\bar{z_1 },d\bar{z_2 }$ then descend to $M$. If $A=(a_{jk})$
 is a $2\times 2$ (constant) Hermitian non-degenerate matrix, then 
$\omega = \Sigma a_{jk}dz_j \wedge d\bar{z_k}$ defines a closed, real, 
(1,1)-form on $M$. So $\omega $ is the
\k form of a \k metric $g$. Moreover, this pseudo-metric is flat. If we choose
$A$ to be indefinite, then $g$ is an indefinite \ke metric on $M$ with
Einstein constant 0.

\ \

{\em Minimal Ruled Surfaces}: Let $S$ be a Riemann surface of genus
$g\ge 2$. Choose a Riemannian metric $h_1$ on $S$ with constant scalar
curvature -1. Fixing an orientation on $S$ there is an almost complex
structure on $S$ (giving the orientation) for which $h_1$ is Hermitian. 
Then $h_1$ is a \ke metric
on $S$ with Einstein constant -1. In the same way construct a \ke metric
$h_2$ on ${\bf CP}^1 $ with Einstein constant 1. Then $h_2 -h_1 $ is a
well defined indefinite \ke metric on $M={\bf CP}^1 \times S$ with
Einstein constant 1.

The general  ruled surface is of the form ${\bf P}(E)$, where
$E$ is a 2-dimensional complex vector bundle over a Riemann surface $S$.
We will later construct indefinite \ke metrics on `most' of these twisted
products (assuming always that the genus of $S$ is greater than 1).

\ \

Now we can state the main result of this paper:

\begin{theorem}Let $M$ be a compact complex surface. If $M$ admits  an
indefinite K\"{a}hler-Einstein metric, then $M$ is one of the following:

 a)  a Complex Torus ; 

 b)  a Hyperelliptic surface ; 

 c)  a Primary Kodaira surface ; 

 d)  a minimal  ruled surface over a curve of genus $g \geq 2$ ; or 

 e)  a minimal surface of class $VII_0 $  with no global spherical shell, and
     
with second Betti number even and positive.

\end{theorem}

{\em Remarks}: No surface of type (e) is known, and it has been
conjectured that they simply do not exist (see \cite{Nakamura}, section 5).
We will display indefinite \ke  metrics with Einstein constant 0 
on the surfaces
(a) , (b)  and (c) (we have already done it for (a)). 
We will also display indefinite \ke  metrics 
with Einstein constant
1 on `most' surfaces in (d); but it is not known if every surface described
by (d) admits such a metric.

\section{Indefinite K\"{a}hler Metrics}

In this section we will find obstructions to the existence of indefinite
\k metrics on a compact complex surface $M$. 

The first thing to note is that
the \k form $\omega $ of such a metric is a symplectic form on $M$.
Since the metric is  indefinite Hermitian  we can find an orthogonal
basis of $T_p M$ of the form $\{ x,Jx,y,Jy\}$ such that $\omega (x,Jx)>0$
and $\omega (y,Jy)<0$. Hence $\omega \wedge \omega $ defines the non-standard 
orientation of $M$.

Moreover, given our symplectic form $\omega $, there exists a {\em compatible
positive almost complex structure} $J_{\omega}$ on $M$ (see for instance 
\cite{Vaisman}, pages 40, 56). Such $J_{\omega}$ is then an almost 
complex structure giving the non-standard orientation of $M$.

{\it Notation}: As usual $b_k (M)$ will denote the $k$-th Betti number 
of the 4-dimensional manifold 
$M$ and $b^+ (M)$ $(b^- (M))$ will denote the dimension of a
maximal subspace of $H^2 (M^4 ,{\bf R})$ where the intersection
form is positive (negative) definite. So $b^+ + b^- =b_2 $ and 
$b^+ -b^- =\tau$, the signature of $M$.

The Todd genus $Todd(M)=1/2(1-b_1 +b^+ )$ of an almost complex manifold $M$ (of
real dimension 4)
is an integer. So $b_1 (M) -b^+ (M)$ is  odd. Assume that $M$ also admits
an almost complex structure compatible with the opposite orientation; 
then $b^+ (M)- b^- (M)=\tau $ is even.
We have proved

\begin{proposition}Let $M$ be a compact complex surface that admits an 
indefinite
\k metric. Then there is an almost complex structure giving the 
non-standard orientation of $M$. In particular, $\tau (M)$ is even.
\end{proposition}

The main obstructions to the existence of indefinite \k metrics on
compact complex surfaces will be obtained using Seiberg-Witten
invariants. These were introduced very recently (see \cite{Witten})
and we will now outline (very roughly) some of the main facts about them
(see \cite{Kronheimer}, \cite{LeBrun}, \cite{Witten} for
details and more general statements). We will follow \cite{LeBrun}.

Let $X$ be a 4-dimensional smooth compact oriented manifold. 
Assume that $X$ admits
an almost complex structure $J$ (compatible with the given orientation)
and that $b^+ (X)\geq 2$. These last conditions are not necessary for
the definition of Seiberg-Witten invariants, but will simplify the description
and we will not need more general results in this work.

Fix an homotopy class $c$ of almost
complex structures on $X$. The Seiberg-Witten invariants will depend only
on the class $c$ on $X$, but Riemannian metrics will be involved in the
construction.

An almost complex structure $J$ provides $TX$ with the structure of a
(2-dimensional) complex vector bundle. Moreover, homotopic almost
complex structures produce isomorphic complex vector bundles. Hence
the homotopy class $c$ provides $TX$ with a canonical complex
structure. Let $T^{1,0} $ be this complex vector bundle and $L={\Lambda}^2
T^{1,0}$. The first Chern class of $L$ is an integral lifting of
$w_2$ (the second Stiefel-Whitney class of $TX$). Then $L$ is a $Spin^c $
structure on $X$, that depends only on the class $c$.

Not every $Spin^c $ structure on $X$ is obtained in this way, and 
Seiberg-Witten invariants can be defined for any $Spin^c$ structure.
But considering only those coming from almost complex structures
will keep things a little simpler.

Any $Spin^c $ structure produces complex 2-dimensional vector bundles
$V^+ $, $V^- $ over $X$ such that $V^+ \bigotimes V^{-*} \cong TX\bigotimes 
{\bf C}$. And choosing a connection $A$ on the line bundle $det (V^+ )$,
it is induced a Dirac operator  

$$D_A : C^{\infty}(V^+ )\rightarrow C^{\infty}(V^-)$$

If the $Spin^c$ structure is given by the line bundle $L$ 
coming from the homotopy class $c$,
 the vector bundles $V^{\pm}$ can be described more
explicitly:

$$\begin{array}{lll}
   V^+ & \cong & \underline{\bf C} \  \  \bigoplus  \  \ L \\ 
   V^- & \cong & T^{1,0} 
\end{array}$$
\noindent where $\underline{\bf C}$ is the trivial line bundle.
So $det(V^+ )=L$.

Given a Riemannian metric $g$ on $X$, we can choose $J$ in $c$
such that $g$ is Hermitian with respect to $J$. In this case the metric $g$
induces a Hermitian metric on $L$. 

The metric $g$ and the orientation of $X$ produce the Hodge operator
$\ast :{\Lambda}^2 T^{\ast}X \rightarrow {\Lambda}^2 T^{\ast}X$; ${\ast}^2=1$
and considering the $\pm 1$ eigenspaces of $\ast$ we have a splitting
${\Lambda}^2 T^{\ast}X = {\Lambda}^+ \oplus {\Lambda}^-$. The sections of
${\Lambda}^+$ are called self-dual 2-forms and the sections of ${\Lambda}^-$
anti self-dual.

Fixing a real self-dual 2-form $\varepsilon \in C^{\infty}({\Lambda}^+ )$ the
(perturbed) Seiberg-Witten equations are

$$ \left\{ \begin{array}{rccc}
                 D_A \Phi  & = & 0 & (1) \\
              iF_A^+ + \sigma (\Phi ) & = & \varepsilon & (2)
           \end{array}\right.  $$

The unknowns are $\Phi \in C^{\infty}(V^+ )$ and $A$, a unitary connection on
$L$. ${F}_{A}^{+}$ is the self-dual part of the curvature of $A$ and
$\sigma :V^+ \rightarrow {\Lambda}^{+}$ is given by

$$\sigma (f,\psi )=(|f|^2 -|\psi |^2 )\omega /4 + Im(\bar{f}\psi )$$
\noindent where $\omega $ is the \k form of $g$.

The Seiberg-Witten equations depend on the pair $(g,\varepsilon )$; the
next task is to extract from them invariants independent of $(g,\varepsilon )$.

The `gauge group' of $C^{\infty}$ maps $u:M \rightarrow S^1 \subset 
{\bf C}$ acts on the space of solutions of (1) by $(A,\Phi )\mapsto 
(A+2d\log u,u\Phi )$. Let ${\cal M}(g,c)$ denote the
space of solutions of (1) modulo this action. We can consider (2) as
a map $\rho : {\cal M}(g,c)\rightarrow C^{\infty}({\Lambda }^+)$.

Let ${\varepsilon }_H $ denote the harmonic part of $\varepsilon \in C^{\infty}
({\Lambda }^+ )$ and $c_1^+$ be the image of $c_1 (L)$ under the orthogonal
projection of $H^2 (M,{\bf R})$ over $H^+ (g)=\{[\varphi ]:\varphi$
is a harmonic
self-dual 2-form (with respect to  $g)\}$.

The pair $(g,\varepsilon )$ is called $excellent$ if
$2\pi c_1^+ \neq [{\varepsilon}_H ]$ and $\varepsilon$ is a regular value
of $\rho $.

Under the assumptions that $b^+ (X)\geq 2$ and that the $Spin^c $ structure
comes from an almost complex structure, we have the following facts:

\ \

a) There exist excellent pairs $(g,\varepsilon )$ and the set of excellent
pairs is path connected.

\ \

b) For any excellent pair $(g,\varepsilon )$, the space of solutions
of the Seiberg-Witten equations modulo de action of the gauge group
is a finite set of points, with a `canonical smooth 
structure'.

\ \

c) Given two excellent pairs the spaces of solutions of the corresponding
Seiberg-Witten equations (modulo de action of the gauge group) are
cobordant.

\ \

Given an excellent pair $(g,\varepsilon )$ let
 $$n_c (N,g,\varepsilon )=\# \{ \mbox{ gauge classes of solutions of (1) and 
(2)}\}\ \mbox{ (mod 2)}$$

We can now define the simplest version of the Seiberg-Witten invariants:

\begin{definition}Let $X$ be a smooth, compact, oriented
 4-dimensional manifold which admits
almost complex structures and such that $b^+ (X)\geq 2$. The (mod 2)
Seiberg-Witten invariant
$n_c (X)$ of $X$ with respect to 
the homotopy class $c$ of almost complex structures
is defined to be $n_c (X,g,\varepsilon )$; where $(g,\varepsilon )$ 
is an excellent
pair.
\end{definition}

\ \

The next theorem of Taubes \cite{Taubes} will give the strongest
tool to show that some compact complex manifolds do not admit
indefinite \k metrics.

\begin{theorem}[Taubes]Let $X$ be a compact, oriented, 4-dimensional manifold
with $b^+ \geq 2$. Let $\omega $ be a symplectic form on $X$ with
$\omega \wedge \omega$ giving the orientation. Then the associated 
homotopy class of almost complex structures $c$ has (mod 2) Seiberg-Witten
invariant 1.
\end{theorem}

This shows that if $M$ is a compact complex surface with $b^- (M)
\geq 2$ and if $M$ admits an indefinite \k metric, then the Seiberg-Witten
invariant of the induced class of almost complex structures is 1.

The following are `classical' results in the theory of
Seiberg-Witten invariants.

\begin{proposition}[\cite{Kronheimer}]Let $g$ be a Riemannian metric on
the smooth compact oriented 4-manifold $X$. Consider the Seiberg-Witten
equations for the pair $(g,0)$. Any solution $(A,\Phi )$ satisfies the
$C^0$ bound
$$\|\Phi \|^2 \leq max(0,-s)$$
\noindent at the points where $\|\Phi \| $ is maximum. Here $s$ is the 
scalar curvature of the Levi-Civita connection of $g$.
\end{proposition}

\begin{theorem} Let $X$ and $Y$ be smooth compact oriented 4-manifolds.
Assume that $b^+ (X)\geq 1$ and $b^+ (Y)\geq 1$. If $c$ is any homotopy
class of almost complex structures on the connected sum
of $X$  and $Y$ (compatible with the given orientation),
then $n_c (X\# Y)=0$
\end{theorem}

Now we can prove

\begin{proposition} Let $M$ be a compact complex surface. If $M$ is 
obtained by blowing up another surface $N$ at one
point and $b^{-}(N)$ is positive, then $M$ does not admit an 
indefinite K\"{a}hler metric. In particular if
$M$ is obtained by blowing up another surface twice, then $M$ does not
admit an indefinite K\"{a}hler metric.

\end{proposition}
Proof: If $M$ is the blown up of $N$, then $M$ is diffeomorphic to the 
connected sum of $N$ with $\overline{{\bf CP}^2}$ (${\bf CP}^2 $ provided
with the non-standard orientation). Reversing the orientations and 
applying the last theorem we get that the Seiberg-Witten
 invariant of $M$ with respect to any
almost complex structure giving the non-standard orientation is 0. If 
$M$ admits an indefinite \k metric, Theorem 2 tells us that for at least
one of those almost complex structures the Seiberg-Witten invariant is 1. 
Hence there is no such a metric on $M$.

\QED 

Now we can study which compact complex surfaces admit indefinite \k
metrics. Recall that a complex surface $M$ is called {\em minimal}
if it does not contain an embedded ${\bf CP}^1$ with self intersection
-1 (i.e. $M$ can not be blown-down). 

\begin{lemma} Assume that the compact complex surface
 $M$ admits an indefinite \k metric. Then $M$ is minimal or is a one-point 
blow-up of either ${\bf CP}^2$ or a fake ${\bf CP}^2$.
\end{lemma} 

{\em Remark}:  A fake ${\bf CP}^2$ is a compact complex surface $M$ such
that $b_1 (M)=0$ and $b_2 (M)=1$. It is known that there are only
finitely many of these surfaces (see \cite[p.136]{Barth}).

\ \

Proof: Assume that $M$ admits an indefinite \k metric. Suppose that $M$
is the blow-up of a surface $N$. Proposition 3 implies that $N$ has
positive definite intersection form (i.e. $b^- (N)=0$) and is minimal. 

If $b_1 (N)=b_1 (M)$ is odd Proposition 1 implies that $M$ does not admit an
indefinite \k metric. So $N$ is of \k type.

Now it is time to check the classification of compact complex surfaces
(see \cite{Barth}, for example). If $Kod(N)$ is 0 or 1, then $c_1^2 (N)=0$
and the  formula $c_1^2 + 8q + b^- =10p_g + 9$ would imply $b^- (N)>0$.
If $Kod(N)=-\infty $ the only possibility is $N={\bf CP}^2$. If
$Kod(N)=2$, then $0<c_1^2 =2c_2 +3\tau \leq 3c_2 $. 
This implies that $b_1 (N)=0$ and $b_2 (N)=1$.

\QED

\begin{lemma} If $b^- (M)=0$ then $M$ admits no indefinite \k metric
(so this is the case for ${\bf CP}^2$, fake ${\bf CP}^2$'s, 
secondary Kodaira surfaces and surfaces of class $VII $ with vanishing second
Betti number).
\end{lemma}

Proof: This is simply because if $\omega $ is the \k form of an
indefinite \k metric on $M$, then $\omega $ defines an element in
$H^2 _{de Rham} (M)$ with ${[\omega ]}^2 <0$.

\QED

\begin{lemma}No indefinite \k metric exists on any $K3$ surface or Enriques
surface.
\end{lemma}

Proof: Enriques surfaces are quotients of $K3$ surfaces. If $g$ is an
indefinite \k metric on an Enriques surface, its pull-back would define
an indefinite \k metric on a $K3$ surface. Hence it is enough to prove
that $K3$ surfaces do not admit indefinite \k metrics.

Let $\overline{M}$ be a $K3$ surface endowed with the non-standard
orientation. Then ${c_1}^2 (\overline{M})=96$ and $b^+ (\overline{M})=19$.
The proof given by S.T. Yau to Calabi's conjecture \cite{yau} shows that 
$M$ admits a scalar flat metric $\hat{g}$. 
Since $c_1^2 (\overline{M})>0$,
we have $c_1^+ (\overline{M})\neq 0$ (check the discussion on Seiberg-Witten
invariants above). This condition assures that no pair of the form $(A,0)$
could be a solution of the Seiberg-Witten equations if $\varepsilon =0$.
And using Proposition 2 we get that for the pair $(\hat{g},0)$ the 
Seiberg-Witten equations have no solution at all.
Hence the pair $(\hat{g},0)$ is excellent 
and the Seiberg-Witten invariants of 
$\overline{M}$ vanish. Then Theorem 2 implies  that $M$ does not admit
indefinite \k metrics.

\QED

\begin{lemma} If $M$ is a surface of class $VII_0 $ with a global spherical
shell (see \cite{Nakamura}) and $b_2 (M)=b^- (M)>0$, then $M$ 
does not admit an indefinite \k
metric.
\end{lemma}

Proof: Such a surface is diffeomorphic to the connected sum
of $S^1 \times S^3 $ with $b_2 (M) $ copies of ${\overline{{\bf CP}}}^2$
(see \cite{Nakamura}). The lemma follows from Theorems 2 and 3.

\QED

By the classification of compact complex surfaces (see \cite{Barth}),
the previous lemmas prove:

\begin{theorem}Suppose $M$ admits an indefinite K\"{a}hler metric. Then

 i)  If $Kod(M)=-\infty $, then $M$ is a ruled surface or is as in 

Theorem 1 (e). 

 ii) If $Kod(M)=0$, then $M$ is a torus, an Hyperelliptic surface or a 

Primary  Kodaira surface.

 iii) If $Kod(M)=1$, then $M$ is minimal. 

 iv) If $Kod(M)=2$, then $M$ is minimal or the blow-up of a fake 
     
${\bf CP}^2$.

\end{theorem}

{\em Remark 1} : Every ruled surface $M$ is of the form ${\bf P}(E)$; where
$\pi : E\rightarrow S$ is a 2-dimensional holomorphic vector bundle
over a Riemann surface $S$. Given a Hermitian metric on $E$ and a
K\"{a}hler form ${\omega }_0 $ on $S$, a sign variation on a well
known form gives $\omega = {\pi }^*({\omega }_0 )- is\partial
\bar{\partial } \log ||W|| $; which, for small $s$, is the K\"{a}hler
form of an indefinite \k metric on $M$.

{\em Remark 2} : The surfaces listed in (ii) do admit indefinite \k
metrics. We will show later that they actually admit indefinite \ke
metrics.

{\em Remark 3} : It is not known which surfaces like (iii) and (iv) admit
indefinite \k metrics . But the product of two Riemann surfaces of genus
$g\geq 2$ belongs to (iv) and the product of an elliptic curve and a
curve of genus $g\geq 2$ belongs to (iii) and both admit an indefinite \k
metric.

\section{Indefinite K\"{a}hler-Einstein Metrics }

We will first compute the Kodaira number of a compact complex surface
that admits an indefinite \ke metric.

\begin{proposition} If $M$ admits an indefinite K\"{a}hler-Einstein
 metric with Einstein
constant $\neq 0$, then $Kod (M)=-\infty $ and $c_{1}^2 <0$.
\end{proposition}

Proof: If $M$ admits such a metric then its Ricci form $\rho =k \omega $ is
everywhere non-degenerate and indefinite. If $\gamma \in {\cal O}
(K^m )$ is not trivial then 
$$\rho = \frac{1}{im} \partial \bar{\partial} \log |\gamma |^2 $$

\noindent would be semi-negative where $|\gamma |$ attains its
maximum. Hence, for all $m>0$, $K^m $ has no non-trivial global
section; and $Kod(M)=-\infty $.

The second assertion follows from the facts that $[\rho ]=2\pi c_1$ and
$\omega \wedge \omega <0$.

\QED

\begin{corollary}If $M$ admits an indefinite \ke metric with Einstein
$\neq 0$, then $M$ is as in (d) or (e) of Theorem 1.
\end{corollary}

\begin{proposition} If $M$ admits an indefinite K\"{a}hler-Einstein
 metric with Einstein
constant 0, then $Kod(M)=0$ and $c_1 (M,{\bf R})=0$.
\end{proposition}

Proof: Suppose that $M$ admits such a metric $g$. Then  
$c_1 (M,{\bf R})=0$ and
 $M$ must be minimal. The only surfaces with Kodaira number $-\infty $ and 
vanishing real first Chern class are the minimal surfaces of class
$VII $  with 0 second
Betti number; which do not admit indefinite \k metrics.
So we can assume that there exists $m>0$ and $\gamma \in {\cal O}(K^m_M )$
non-trivial.
Let $\widetilde{M}  $ be the universal covering
of $M$. The pull-back of $g$ gives an indefinite \ke  metric on 
$\widetilde{M}$ 
(with Einstein
constant 0). Since this metric is Ricci flat, there are holomorphic 2-forms
of constant length in a neighborhood of any point (this fact is usually
stated only in the Riemannian  case, but it  is not difficult to check
that the proof also works in the indefinite case). Since 
$\widetilde{M}$ is simply
connected it then admits a global non trivial holomorphic 2-form  $\varphi $ of
constant length. The pull back $\hat{\gamma }$ of $\gamma $ can be written
$\hat{\gamma }=f{\varphi }^m $ for some holomorphic function $f$ on 
$\widetilde{M}$
of bounded length. Hence $f$ is constant and $\| \gamma \| $ is constant.
Then $\gamma $ is never zero and $K^m_M $ is trivial. It follows that
$Kod(M)=0$.

\QED

\begin{corollary}If $M$ admits an indefinite \ke metric with Einstein
constant 0, then $M$ is as is (a), (b) or (c) of Theorem 1.
\end{corollary}

By now we have already proved Theorem 1. The only thing remaining is
to construct the promised examples of indefinite K\"{a}hler-Einstein metrics.

\ \

{\em Hyperelliptic surfaces} : 
It is shown in \cite[p.585]{Griffiths} that any hyperelliptic surface $M$ is
of the form $M=F\times C/G $, where $F$ and $C$ are elliptic curves and $G$ is
finite group  of fixed-point-free automorphisms of $F\times C$. Moreover,
let $F={\bf C}/\Lambda$ with $\Lambda =<1,\tau >$; 
then $G=< \phi , \varphi >$, where
$\phi $ is of the form 
$\phi (z,w)=(z+ \tau /m,e^{2k\pi i/m})$ and
$\varphi $ is a translation of order $m$.
If $z,w$ are the standard
holomorphic coordinates in ${\bf C}^2$ then $dz\wedge d\bar{z} -
dw\wedge d\bar{w}$ is the \k form of an indefinite \k flat metric 
(on ${\bf C}^2$). This form is invariant through translations and so
projects to a $(1,1)$-form on $F\times C$.
A direct computation shows that this form is invariant through $\phi $ and
$\varphi $ and hence defines a $(1,1)$ on $M$; this is the \k form of
an indefinite \ke metric on $M$ with Einstein constant 0.

\ \

{\em Primary Kodaira surfaces} : 
As described in \cite[p.786]{Kodaira} such a surface $M$ is of
the from $M={{\bf C}}^2 /G$ where $G=<g_1 ,g_2 ,g_3 ,g_4 >$,
each $g_i $ is an affine automorphism of ${{\bf C}}^2 $ and $G$ is
fixed point free. More precisely each $g_i $ is of the form

$$g_i (w_1 ,w_2 )=(w_1 +{\alpha }_i ,w_2 + {\bar{\alpha }}_i w_1 + 
{\beta}_i ) \ \ ,\ \  {\alpha }_i , {\beta }_i \in {\bf C}$$

Consider the (1,1)-form $\gamma =-(w_1 + {\bar{w}}_1 )dw_1 \wedge d{\bar{w}}_1
 + dw_1 \wedge
d{\bar{w}}_2 + dw_2 \wedge d{\bar{w}}_1 $ on ${\bf C}^2$.
Direct computations show that
$\gamma $ defines an indefinite \k flat metric on ${\bf C}^2 $  
and is invariant
under the  $g_i $'s. Hence $\gamma $ induces an indefinite \ke  metric on $M$
with Einstein constant 0.

\ \

{\em Minimal irrational ruled surfaces} : Now let $M={\bf P}(E) $, where
$E$ is a 2-dimensional holomorphic vector bundle over a curve $S$ of
genus $g\geq 2$. We will construct indefinite \ke metrics on
$M$ when the bundle $E$ is stable or the direct sum of two line bundles
of the same degree (see \cite{Kobayashi}, \cite{Nara}). 

Note that given vector bundles $E$ and  $\widehat{E}$,
 ${\bf P}(E)$ and ${\bf P}(\widehat{E})$ are isomorphic if and only if
$\widehat{E}=E\otimes L$ for a line bundle $L$;  and that $\widehat{E}$
verifies any of the conditions above  if and only if $E$ does. 
So both conditions are really
properties of $M$.

Consider $M$ as a ${\bf CP}^1$-bundle over $S$. Let $(U_i )^N_{i=1}$ be an
open cover of $S$ and $g_{ij}:U_i \bigcap U_j \rightarrow Gl(2,{\bf C})$ 
be a set of 
transition functions for $E$. 
Then $[g_{ij}]~:~U_i ~\bigcap ~U_j ~\rightarrow ~{\bf P}~Gl(2,{\bf C})$ are
transition functions for $M$. Under the conditions stated above, 
Narasimhan and Seshadri \cite{Nara} proved
that $M$ admits constant transition functions in ${\bf P}(U2)$.  
Let $g_1$ be the Fubini-Study metric on ${\bf CP}^1 $; then 
$g_1$ is a \ke metric  on ${\bf CP}^1$, invariant through the action
of ${\bf P}(U2)$. Renormalize $g_1$ so that the Einstein constant is 1
and let  $g_2$ be a K\"{a}hler-Einstein metric on $S$ with 
Einstein constant -1. Then $g_1 -g_2 $ is invariant through the transition
functions and so defines an  indefinite K\"{a}hler-Einstein metric on
$M$ with Einstein constant 1. 

\ \

{\em Remark:} In \cite[p.395]{ramanan} M.S. Narasimhan and S. Ramanan proved
that every vector bundle (over a curve of genus greater than 1) can be
`approximated' by stable vector bundles. A little more precisely, every
vector bundle is contained in an analytic family of vector bundles for
which the set of stable bundles is open and dense. 

The cases considered above therefore contain `most' of the minimal ruled
surfaces (over curves of genus greater than 1).

\ \
 
\noindent
{\bf Acknowledgements:} The author most gratefully thanks Claude LeBrun
for his continuing help and support during the preparation of this paper.

\end{document}